\documentclass[a4paper,12pt]{article} 
\newdimen\paperhight
\usepackage{amsmath}
\usepackage{amssymb}
\usepackage{amsfonts}
\usepackage[usenames]{color}

\newcommand{\dsp}{\displaystyle}

\newcommand{\hf}{\frac{1}{2}}

\newcommand{\pr}{\par \vspace{3mm}\noindent [{\bf Proof}] \qquad}
\newcommand{\prend}{\hfill \qed \par \vspace{3mm}}

\newcommand{\qed}{\quad\hbox{\rule[-2pt]{3pt}{6pt}}\par\vspace{3mm}}

\newcommand{\1}{{\bf 1}} 
 
\newcommand{\C}{\mathbb C} 
\newcommand{\Z}{\mathbb Z}

\newcommand{\N}{\mathbb N}

\newcommand{\CG}{{\cal G}}
\newcommand{\CH}{{\cal H}}
\newcommand{\CJ}{{\cal J}}

\newcommand{\CL}{{\cal L}}
\newcommand{\CP}{{\cal P}}

\newcommand{\FS}{{\cal J}}
\newcommand{\al}{\alpha}
\newcommand{\be}{\beta}

\newcommand{\End}{{\rm End}}

\newtheorem{thm}{Theorem}
\newtheorem{prn}[thm]{Proposition}

\newtheorem{cry}[thm]{Corollary}
\newtheorem{lmm}[thm]{Lemma}
\newtheorem{rmk}[thm]{Remark}
\newcommand{\op}[1]{L^{ij}_{r}(#1)}
\newcommand{\ome}{\omega^{ij}_{r}}
\newcommand{\bs}[4]{:v^{(#1)}{(#2)}v^{(#3)}{(#4)}:_{r} }

\newcommand{\omm}[2]{\omega^{#1#2}_{r}}
\newcommand{\del}[2]{\delta_{#1,#2}}
\begin{document}

\title{Deformation of central charges, 
vertex operator algebras whose Griess algebras are Jordan algebras}
\author{Takahiro Ashihara\\
Graduate School of Pure and Applied Sciences,\\
University of Tsukuba, \cr
Tsukuba, 305 Japan.\cr
megane2@math.tsukuba.ac.jp \\
\\
Masahiko Miyamoto\\
Institute of Mathematics, University of Tsukuba, \cr
Tsukuba, 305 Japan.\cr
miyamoto@math.tsukuba.ac.jp \\}
\date{}
\maketitle
\begin{abstract}
If a vertex operator algebra $V=\oplus_{n=0}^{\infty}V_n$ 
satisfies $\dim V_0=1, V_1=0$, then $V_2$ has a commutative (nonassociative) algebra 
structure called Griess algebra. One of the typical examples of commutative 
(nonassociative) algebras is a Jordan algebra. For example, 
the set $Sym_d(\C)$ of symmetric matrices of degree $d$ becomes a 
Jordan algebra. On the other hand, 
in the theory of vertex operator algebras, central charges 
influence the properties of vertex operator algebras. In this paper, 
we construct vertex operator algebras with central charge $c$ and its 
Griess algebra is isomorphic to $Sym_d(\C)$ for any complex number $c$.\\
\\
Keyword: Vertex operator algebra; Jordan algebra; central charge.
\end{abstract}

\section{Introduction}
A vertex operator algebra (shortly VOA) is a quadruple $(V,Y( ,z),\1,\omega)$ 
consisting of a $\N$-graded vector space $V=\oplus_{n=0}^{\infty}V_n$, a set of 
vertex operators $Y(v,z)=\sum_{n\in \Z}v_nz^{-n-1}\in \End(V)[[z,z^{-1}]]$ 
for $v\in V$, and two special elements, 
a vacuum $\1$ and a Virasoro element $\omega$, and they satisfy several conditions 
called locality and associativity, etc. 
A vertex operator algebra is originally defined as a rigorous mathematical 
object corresponding to (a chiral algebra of) a conformal field theory. 
However it also has a rich algebraic structure. For example, 
if a VOA $V=\oplus_{n=0}^{\infty}V_n$
satisfies $\dim V_0=1$, then the space $V_1$ of weight $1$ becomes 
a Lie algebra by $0$-product $v_0u$ for $v,u\in V_0$, and 
if $\dim V_0=1$ and $V_1=0$, then the space $V_2$ of weight $2$ 
becomes a (nonassociative) commutative algebra (called Griess algebra) by 
$1$-product $v_1u\in V_2$ for $v,u\in V_2$, where 
$Y(v,z)=\sum v_nz^{-n-1}\in \End(V)[[z,z^{-1}]]$ is a vertex operator of $v$. 
The most famous Griess algebra is the monstrous Griess algebra of the moonshine 
vertex operator algebra $V^{\natural}$, whose full 
automorphism group is the monster simple group.  

On the other hand, one of the most typical (nonassociative) commutative algebras 
is a Jordan algebra. For example, 
${\rm Sym}_d(\C)$, the set of symmetric matrices of degree $d$, becomes a 
Jordan algebra by a product $A\ast B=\frac{1}{2}(AB+BA)$. 
So, it is a natural question if there is a VOA $V$ whose Griess algebra 
is isomorphic to a given Jordan algebra $\CG$ or not. 

As an answer to this question, Lam \cite{Lam} constructed such VOAs as 
subVOAs of $M(1)^{\otimes d}$ for the most of complex simple Jordan algebras, where 
$M(1)$ is a VOA associated with a Heisenberg algebra and 
$M(1)^{\otimes n}$ denotes a tensor product of $n$ copies of $M(1)$, see \cite{FLM}. 

Recently, Jordan algebras win another attention because of 
connection with symmetric cones and Zeta functional equations \cite{FK}. 
In this case, an inner product of Jordan algebra plays an important role. 
For a VOA $V=\oplus_{n=0}^{\infty} V_n$ with $\dim V_0=1$ and $V_1=0$, 
a Griess algebra $V_2$ also has a natural 
invariant bilinear form $\langle , \rangle$ defined by 
$\langle a,b\rangle\1=a_3b$ for $a,b\in V_2$ and $2\langle \omega,\omega\rangle$ 
is called a central charge of $V$, where $\omega$ is a Virasoro element of $V$ and 
$\1$ is the vacuum element. 
Moreover, in the theory of vertex operator algebra, 
many properties of VOAs are influenced by central charges. For example, it is well known that 
if all modules are completely reducible, then central charges are all rational numbers. 
Therefore, it fertilizes the theory of vertex operator algebras 
to construct many examples with arbitrary central charges. 

For example, there are infinitely many VOAs $V_{\CG}$ having 
same Lie algebra as the weight one space $V_1$, but with different 
central charges \cite{FZ}. 
These VOAs have close connection with 
Lie algebras. On the contrary, a VOA without weight one space 
keeps Lie algebra in the background and 
has many interesting relations with finite group and other finite dimensional algebras. 

The purpose of this paper is to fuse these sights. 
Since the central charges of VOAs constructed by Lam are 
the sizes of Jordan algebras, that is, natural numbers, 
we would like to deform central charges and 
show a new construction of vertex operator algebras, 
whose Griess algebras are Jordan algebras 
and having any complex number as a central charge. 
Namely, we will prove the following theorem. 

\begin{thm}
For any complex number $r$ and a natural number $d$, 
there is a VOA $V$ with central charge $rd$ such that 
its Griess algebra is isomorphic to a Jordan algebra ${\rm Sym}_d(\C)$.
\end{thm}

\begin{rmk}  If $d=1$, then ${\rm Sym}(1,\C)\cong \C$ and our construction 
coincides with 
a construction of Virasoro VOAs with central charge $r$. Therefore, one may consider 
our VOAs as natural extensions of Virasoro VOAs. 
\end{rmk}

\section{Deformation of central charges}

\subsection{Deformed Lie algebra}
Let $H$ be a $d$-dimensional vector space and $r\in \C$. Viewing it as a 
commutative Lie algebra with an invariant bilinear form, 
we can define its affinization:
$$ \begin{array}{l}
\hat{H}=H[t,t^{-1}]\oplus \C{\bf c}, \cr
\end{array}$$
where Lie products are given by
$$ \begin{array}{l}
[h_1t^m, h_2t^n]=\delta_{m+n,0}m<h_1,h_2>{\bf c} \mbox{ for }h_1,h_2\in H \cr
[{\bf c},\hat{H}]=0 
\end{array}$$
and consider the universal enveloping algebra $U(\hat{H})$ of $\hat{H}$ and 
its factor ring: 
$$U(\hat{H})/({\bf c}-r)=\C \oplus H[t,t^{-1}]\oplus S^2(H[t,t^{-1}])\oplus ... $$
divided by an ideal generated by ${\bf c}-r$, 
where $S^n(U)$ is the weight $n$ space of a symmetric tensor algebra 
$S(U)=\C\oplus U\oplus S^2(U)\oplus \cdots$ of $U$. 
Viewing this algebra as a Lie algebra with a Lie product 
$[\al,\be]=\al\be-\be\al$, it is easy to see that 
$$\CL_r=S^2(H[t,t^{-1}])\oplus \C   \quad(\subseteq U(\hat{H})/({\bf c}-r))$$ 
is a subalgebra. We sometimes denote subspaces 
$S^2(H[t,t^{-1}])$ and $\C$ of $\CL_r$ by $\CL_r'$ and $\CL_r"$, respectively. 
We introduce a new Lie product in $\CL_r$ by 
$$ [\al,\be]_r=\frac{1}{r}[\al,\be].$$
The case $r=1$ is standard, that is, when we construct a VOA of free bozon, 
$U(\hat{H})/({\bf c}-1)$ appears.  
We note that $\CL_r$ coincides with $\CL_1$ as a vector space, but has a 
different Lie algebra structure for each $r$. 
An important fact is that there is another way to define this Lie algebra $\CL_r$  
by introducing a new Lie product
$$ [\al,\be]=[\al,\be]_1'+r[\al,\be]_1" \quad \forall \al,\be\in \CL_1, $$
 in $\CL_1=\CL_1'\oplus \CL_1"$, where $[\al,\be]_1'$ and 
$[\al,\be]_1"$ denote projections of $[\al,\be]_{1}$ into $\CL_1'$ and $\CL_1"$, respectively.  \\

We next define a Fock space $M_r$. 
Let $\{v^{(1)},...,v^{(d)}\}$ be an orthonormal basis of $H$.
We denote $v^{(i)}t^m\in H[t,t^{-1}]$ by $v^{(i)}(m)$ for simplicity. 
Clearly, we have:

\begin{lmm}\label{lmm:2} 
$$ \{v^{(i)}(m)v^{(i)}(n)\mid  m\leq n, i=1,...,n\}\cup \{v^{(i)}(m)v^{(j)}(n)\mid 
1\leq i<j\leq d, 
m,n\in \Z\}$$
is a basis of $\CL_r'=S^2(H[t,t^{-1}])$.
\end{lmm}

In order to distinguish a new product $[\cdot,\cdot]_r$ from $\CL_1$, 
we denote the above basis by 
$$:v^{(i)}(m)v^{(j)}(n):_r.$$ 

In $\CL_r$, 
$$\begin{array}{l}
\CL_r^+=\langle :v^{(i)}(m)v^{(j)}(n):_r \mid m\geq 0\mbox{ or }n\geq 0 \rangle
\oplus \C\qquad \mbox{ and }\cr
\CL_r^-=\langle :v^{(i)}(m)v^{(j)}(n):_r \mid m\lneq 0,n\lneq 0\quad \rangle
\end{array}$$ 
are sub Lie algebras of $\CL_r$ and 
$$\CL_r=\CL_r^-\oplus\CL_r^+.$$ 
We note that $\CL_r^-$ is a commutative Lie algebra. By 
the Poincar\'{e}-Birkhoff-Witt theorem, we have:
$$ U(\CL_r)=U(\CL_r^-)\otimes_{\C}U(\CL_r^+),$$ 
where $U(\CH)$ denotes the universal enveloping algebra of $\CH$. 

We next define a $U(\CL_r^+)$-module $\C\1$ by the conditions: \\
(1) $:v^{(i)}(m)v^{(j)}(n):_r\1=0$ if $m\geq 0$ or $n\geq 0$ and \\
(2) the constants $s\in \C$ act as scalar times. \\
Then by inducing, we can define a $U(\CL_r)$-module: 
$$M_r=U(\CL_r)\otimes_{U(\CL_r^+)}\C\1\cong U(\CL_r^-)\1.  $$

Define the grading into $M_r$ by 
$$\begin{array}{l}
{\rm deg}(:v^{(i)}(m)v^{(j)}(n):_r)=-m-n \mbox{  and  }\cr
{\rm deg}(:v^{(i_1)}(m_1)v^{(j_1)}(n_1):_r\cdots :v^{(i_k)}(m_k)v^{(j_k)}(n_k):_r)\1=
-(\sum_{i=1}^k m_i+n_i) 
\end{array}$$
and set 
$$ (M_r)_n=\{v\in M_r \mid {\rm deg}v=n \}.$$  
It is easy to see the following lemma. 

\begin{lmm}\label{lmm:3} 
$$M_r=\oplus_{n=0}^{\infty}(M_r)_n$$ 
and $(M_r)_0=\C{\bf 1}$, 
$(M_r)_1=0$, and $(M_r)_2=\oplus_{i\geq j}\C:v^{(i)}(-1)v^{(j)}(-1):_r\1$. 
\end{lmm}

\subsection{Vertex operators}
Clearly, $\CL_r$ acts on $M_r$ faithfully and we may view 
$\CL'\oplus \CL"\cong \CL'_r\oplus \CL_r"\subseteq 
{\rm End}(M_r)$. The useful fact is that 
the commutators of $\CL_r$ are still in $\CL_r$, which plays an important role. 

In imitation of the Sugawara-construction, 
we define operators $L_r^{ij}(n)$ on $M_r$ as follows:
$$ \begin{array}{l}
L_r^{ii}(n):=\hf \sum_{n-h\leq h}:v^{(i)}(n-h)v^{(i)}(h):_r
+\hf\sum_{h<n-h}:v^{(i)}(h)v^{(i)}(n-h):_r  \quad \mbox{ and } \vspace{2mm} \cr 
L_r^{ij}(n):=\hf \sum_{h\in\mathbb{Z}}:v^{(i)}(n-h)v^{(j)}(h):_r \qquad \mbox{if }i\lneq j.
\end{array}$$ 

We note that all terms in $L_r^{ij}(n)$ are in $\CL'_r$. 
Although each $L_r^{ij}(n)$ is an infinite sum of elements in $\CL'_r$, 
all but finitely many act as zero to every element in $M_r$ and so 
the sums are well-defined as elements in ${\rm End}(M_r)$. 

\begin{prn}
We have the following properties:
$$\begin{array}{l}
[L_r^{ii}(m),L_r^{ii}(n)]=(m-n)L_r^{ii}(m+n)+\delta_{m+n,0}\frac{m^3-m}{12}r, \cr
[L_r^{ij}(m),L_r^{st}(n)]=0\mbox{ if }\{i,j\}\cap \{s,t\}=\emptyset, \cr
[L_r^{ii}(m),L_r^{ij}(n)]=\frac{1}{2}\dsp\sum_{k\in\Z}(m-k)\bs{i}{k}{j}{m+n-k},  \cr
[L_r^{ij}(m),L_r^{ij}(n)]=\frac{m-n}{2}(L_r^{ii}(n+m)+L_r^{jj}(n+m))
+\delta_{n+m,0}\frac{m^3-m}{24}r, \cr
[L_r^{ij}(m),L_r^{jk}(n)]=\frac{1}{4}\sum_{l\in \Z} l:v^{(i)}(m-l)v^{(k)}(n+l):_r  
\mbox{ for } i\neq j\neq  k\neq i, {\rm and}\cr
[\sum_{i=1}^{d}L^{ii}_r(-1),L^{st}_r(m)]=(-1-m)L^{st}_r(-1+m),
\end{array}$$
where $[A,B]=AB-BA$.
\end{prn}

\pr 
As it is well-known, we are able to obtain the above results for the original case $r=1$ 
by the direct calculation. 
From the definition of $L_r^{ij}(n)$, since it is an infinite sum of elements in 
$\CL_r'$, we have:
$$[L_r^{ij}(n), L_r^{st}(m)]=[L_1^{ij}(n),L_1^{st}(m)]'+r[L_1^{ij}(n),L_1^{st}(m)]" $$ 
by viewing $\CL_1=\CL_r\subseteq {\rm End}(M_r)$ as vector spaces, 
where we have to note that $[,]'$ and $[,]"$ in the right side denote  
infinite sums of elements in $\CL_r'$ and $\CL_r"$, respectively. 
Therefore, we have the desired conclusion for any $r$. 
\prend

As a corollary, we have:

\begin{cry}
If we define 
$\omega_r^{ij}(z)=\sum_{n\in \Z} L_r^{ij}(n)z^{-n-2}$ for $i\leq j$, 
then all $\omega_r^{ij}(z)$ satisfy local commutativity with each other 
including itself. 
\end{cry}

\pr
As it is well-known, we have 
$$0=(x-z)^4[\omega_1^{ij}(x),\omega_1^{st}(z)]. $$
Using an expression: 
$$\begin{array}{rl}
[\omega_1^{ij}(x),\omega_1^{st}(z)]&=[\omega_1^{ij}(x),\omega_1^{st}(z)]'
+[\omega_1^{ij}(x),\omega_1^{st}(z)]" \cr
&\in \CL'[[x,x^{-1},z,z^{-1}]]\oplus \CL"[[x,x^{-1},z,z^{-1}]], 
\end{array}$$ 
we have 
$$\begin{array}{rl}
0=&(x-z)^4[\omega_1^{ij}(x),\omega_1^{st}(z)]' \cr
0=&(x-z)^4[\omega_1^{ij}(x),\omega_1^{st}(z)]". 
\end{array}$$
Therefore, we have 
\begin{align*}(x-z)^4[\omega_r^{ij}(x),\omega_r^{st}(z)]=
(x-z)^4[\omega^{ij}(x),\omega^{st}(z)]'
+r(x-z)^4[\omega^{ij}(x),\omega^{st}(z)]"=0.\end{align*} 
\prend
We set
\begin{align*}
\CJ&:=\{\ome~|~1\leq i\leq j\leq d\},\\
V_{\CJ}&:=<L^{i_{1},j_{1}}_{r}(m_{1})\cdots L^{i_{k},j_{k}}_{r}(m_{k}){\1}~|~k\in {\N},1\leq i_{l}\leq j_{l} \leq d,1\leq l\leq k~>,\\
\end{align*}
where $\ome=\op{-2}\1$, and define a map $Y_{0}(\cdot,x)$ from $\CJ$ to End($V_{\CJ}$)[[$x,x^{-1}$]] as follows;
\begin{align*}
Y_{0}(\ome,z):=\ome(z)=\sum_{m\in{\Z}}{\op{m}}x^{-m-2}.
\end{align*}
\begin{lmm}
$Y_{0}$ can be uniquely extended to a linear map $Y$ from $V_{\CJ}$ to End($V_{\CJ}$)[[$x,x^{-1}$]] such that $(V_{\CJ},Y,{\1})$ carries the vertex algebra structure and $L_{r}(-1):=\dsp\sum_{i=1}^{d}L^{ii}_{r}(-1)$ is a operator on $V_{\CJ}$ such that for $v\in V_{\CJ}$
\begin{align*}
[L_{r}(-1),Y(v,x)]=\frac{d}{dx}Y(v,x).
\end{align*}
\end{lmm}
\pr 
By the definition of the operator $\op{m}$, we have
\begin{center}
${\op{m}}{\1}=0$~~~~~if $m\geq -1$.
\end{center}
Hence, we obtain
\begin{align*}
Y_{0}(\ome,x){\1}\in End(V_{\CJ})[[x]]~,~\lim_{z\rightarrow 0}Y_{0}(\ome,x){\1}={\op{-2}}{\1}=\ome
\end{align*}
and
\begin{align*}
L_{r}(-1)\1=0.
\end{align*}
By Proposition 5 we have
\begin{align*}
[L_{r}(-1),Y_{0}(\ome,x)]&=\sum_{m\in \Z}(-m-2)\op{m}x^{-m-3}\\
&=\frac{d}{dx}Y_{0}(\ome,z).
\end{align*}
Therefore, using Theorem 5.7.1 in [8], this lemma holds.
\prend
It is clear by Proposition 5 and Lemma 7 that $V_{\CJ}$ is a conformal vertex algebra with the vacuum vector $\1$, the virasoro element $\omega_{r}:=\dsp\sum_{i=1}^{d}\omega_{r}^{ii}$ , and a central charge $dr$. Since  the element $\bs{i}{m}{j}{n}\in M_{r}$ satisfies the condition
\begin{align*}
&[L_{r}(0),\bs{i}{m}{j}{n}]=(-m-n)\bs{i}{m}{j}{n},
\end{align*}
we obtain the following lemma.
\begin{lmm}
The grading of $M_{r}$ coincides with the $L_{r}(0)$-eigenspace decomposition of $M_{r}$. 
\end{lmm}
\prend
Let $(V_{\CJ})_{m}$ be the homogeneous subspace of weight $m$ of $V_{\CJ}$. By Lemma 8, it is clear that dim$(V_{\CJ})_{m}$ is finite.  
\begin{thm}
$(V_{\CJ},Y,\1,\omega)$ is a vertex operator algebra with the central charge $dr$. In particular,
\begin{align}
&(V_{\CJ})_{n}=\{0\}~~{\rm for~all}~~n\in\mathbb{Z}_{n\lneq 0},\\
&(V_{\CJ})_{0}=\C\1,\\
&(V_{\CJ})_{1}=\{0\},\\
&((V_{\CJ})_{2},1-product)\cong{\rm Sym}_{d}(\C).
\end{align}
\end{thm}
\pr (1)-(3) are easy. We see by the direct calculation that
\begin{align*}
\ome=\op{-2}\1=\frac{1}{2}\bs{i}{-1}{j}{-1}.
\end{align*}
Thus, $\{\ome~|~1\leq i\leq j\leq d~\}$ forms the basis of $(V_{\CJ})_{2}$. In particular, $(V_{\CJ})_{2}=(M_{r})_{2}$. Finally, define a map $g$ from Sym$_{d}(\C)$ to $(V_{\CJ})_{2}$ as follows:
\begin{align*}
g(E^{ij})=\ome,
\end{align*}
where $E^{ij}$ is the matrix of degree $d$ with $i,j$-entry and $j,i$-entry being $1$ and all other entries $0$. Calculating $1$-product and the product introduced in Introduction of Sym$_{d}(\C)$, we have
\begin{align*}
(\omm{i}{j})_{1}(\omm{s}{t})=\frac{1}{2}(\del{i}{s}\omm{j}{t}+\del{i}{t}\omm{j}{s}+\del{j}{s}\omm{i}{t}+\del{j}{t}\omm{i}{s})
\end{align*}
and
\begin{align*}
E^{ij}*E^{st}=\frac{1}{2}(\del{i}{s}E^{jt}+\del{i}{t}E^{js}+\del{j}{s}E^{it}+\del{j}{t}E^{is}).
\end{align*}
Therefor, $g$ is an algebra isomorphism and (4) holds.
\prend
Thus, $V_{\CJ}$ is the desired VOA in Theorem 1 in Introduction.\\

We will denote the above VOA  $V_{\CJ}$ by $V_{\CJ B}(d,dr)$ according to 
the classification of Jordan algebra by Albert \cite{A}, where $dr$ is the 
central charge of $V_{\CJ B}(d,dr)$.\\ 

Since $(V_{\CJ B}(d,dr))_2\cong {\rm Sym}_d(\C)$ is simple, all nontrivial 
ideals of $V_{\CJ B}(d,dr)$ are contained in 
$\oplus_{n\geq 3} (V_{\CJ B}(d,dr))_n$ and 
so $V_{\CJ B}(d,dr)$ has the unique maximal ideal $\CP$. 
In particular, we have a simple VOA $V_{\CJ B}(d,r)/\CP$ with the same Griess 
algebra $(V_{\CJ B}(d,dr)/\CP)_2\cong {\rm Sym}_d(\C)$.

\begin{prn} 
${\rm Aut}(V_{\CJ B}(d,dr))\cong {\rm Aut}(V_{\CJ B}(d,dr)/\CP) 
\cong O_d(\C)/\langle \pm 1\rangle$. \end{prn}

\pr

Clearly, $O_d(\C)$ acts on ${\rm VA}(\FS)$ and also on $M_r$ and 
$\pm 1$ acts trivially on the $M_r$. 
Since $O_d(\C)$ is compatible with $f$, $O_d(\C)$ acts on $V_{\CJ B}(d,dr)$. 
Furthermore, since $O_d(\C)/\langle\pm 1\rangle$ acts 
on ${\rm Sym}_d(\C)$ faithfully, we have 
$O_d(\C)/\langle\pm 1\rangle\subseteq {\rm Aut}(V_{\CJ B}(d,dr)$. 
On the other hand, if there is 
$g\in {\rm Aut}(V_{\CJ B}(d,dr))-O_d(\C)/\langle\pm 1\rangle$, 
then we may assume that $g$ fixes all elements in $(M_r)_2$ since 
$O_d(\C)/\langle \pm 1\rangle={\rm Aut}({\rm Sym}_d(\C))$. 
However, since ${\rm VA}(\FS)$ is generated by operators of 
$(M_r)_2$, we have a contradiction. 
\prend

\end{document}